\documentclass[10pt]{amsart}
\usepackage{amssymb, amsthm, amsmath}
\usepackage{xy}
\usepackage{epsfig}
\usepackage{psfrag}

\psfrag{lambda}{$\lambda$}
\theoremstyle{plain}
\newtheorem{prop}{Proposition}
\newtheorem{thm}{Theorem}
\newtheorem{cor}{Corollary}

\newtheorem{conj}{Conjecture}
\theoremstyle{remark}
\theoremstyle{definition}
\newtheorem{remark}{Remark}

\newtheorem{exa}{Example}

\newcommand{\CC}{{\mathcal M}}
\newcommand{\II}{{\mathcal I}}
\newcommand{\bP}{{\mathbb P}}

\newcommand{\Z}{{\mathbb Z}}

\newcommand{\C}{{\mathbb C}}
\newcommand{\D}{{\mathcal D}}
\newcommand{\cO}{{\mathcal O}}
\newcommand{\X}{{\mathfrak X}}
\newcommand{\Y}{{\mathfrak Y}}
\newcommand{\zZ}{{\mathcal Z}}
\newcommand{\W}{{\mathcal W}}
\newcommand{\s}{\sigma}
\newcommand{\la}{\lambda}
\newcommand{\m}{\mu}
\newcommand{\n}{\nu}
\newcommand{\cc}{\negmedspace:\negmedspace}
\newcommand{\bull}{{\scriptscriptstyle \bullet}}

\DeclareMathOperator{\Span}{Span}
\DeclareMathOperator{\Ker}{Ker}

\newcommand{\ssm}{\smallsetminus}
\newcommand{\gequ}{\geqslant}
\newcommand{\lequ}{\leqslant}

\newcommand{\ra}{\rightarrow}

\newcommand{\ov}{\overline}

\newcommand{\wt}{\widetilde}

\newcommand{\pzp}[1]{\includegraphics[scale=.9]{#1.ps}}
\newcommand{\pzs}{\hspace{4mm}}

\begin{document}

\title[Gromov-Witten invariants on Grassmannians]
{Gromov-Witten invariants on Grassmannians}
\author{Anders Skovsted Buch, Andrew Kresch, and Harry Tamvakis}
\date{April 24, 2003}
\thanks{The authors were supported in part by NSF Grant DMS-0070479
  (Buch), an NSF Postdoctoral Research Fellowship (Kresch),
  and NSF Grant DMS-0296023 (Tamvakis).}
\address{Matematisk Institut, Aarhus Universitet, Ny Munkegade, 8000
  {\AA}rhus C, Denmark}
\email{abuch@imf.au.dk}
\address{Department of Mathematics, University of Pennsylvania,
209 South 33rd Street,
Philadelphia, PA 19104-6395, USA}
\email{kresch@math.upenn.edu}
\address{Department of Mathematics, Brandeis University - MS 050,
P. O. Box 9110, Waltham, MA
02454-9110, USA}
\email{harryt@brandeis.edu}

\begin{abstract} 
We prove that any three-point
genus zero Gromov-Witten invariant 
on a type $A$ Grassmannian is equal to a classical
intersection number on a two-step flag variety. We also give 
symplectic and orthogonal analogues of this result; in these cases
the two-step flag variety is replaced by 
a sub-maximal isotropic Grassmannian. Our theorems are applied, 
in type $A$, to formulate a conjectural quantum Littlewood-Richardson rule,
and in the other classical Lie types, to obtain new proofs of the main
structure theorems for the quantum cohomology of Lagrangian and 
orthogonal Grassmannians.
\end{abstract}

\maketitle

\section{Introduction}
\label{intro}

\noindent
The central theme of this paper is the following result: any three-point
genus zero Gromov-Witten invariant 
on a Grassmannian $X$ is equal to a classical
intersection number on a homogeneous space $Y$ of the same Lie type.
We prove this when $X$ is a type $A$ Grassmannian, and, in types $B$, $C$,
and $D$, when $X$ is the Lagrangian or orthogonal Grassmannian
parametrizing maximal isotropic subspaces in a complex vector space
equipped with a non-degenerate skew-symmetric or symmetric form. 
The space $Y$ depends on $X$ and the degree of the Gromov-Witten 
invariant considered.
For a type $A$ Grassmannian, $Y$ is a two-step flag
variety, and in the other cases, $Y$ is a sub-maximal isotropic
Grassmannian.  

Our key identity for Gromov-Witten invariants is based on an
explicit bijection between the set of rational maps
counted by a Gromov-Witten invariant and the set of points in the intersection 
of three Schubert varieties in the homogeneous space $Y$.  
The proof of this result uses no moduli spaces of maps and 
requires only basic algebraic geometry.
It was observed in \cite{Buch} that the intersection and linear span
of the subspaces corresponding to points on a curve in a Grassmannian,
called the kernel and span of the curve, have dimensions which are
bounded below and above, respectively. 
The aforementioned bijection sends a rational map $\bP^1\ra X$
to the pair consisting of the kernel and span of its image.
(In types $B$, $C$, and $D$, the span is
the orthogonal complement of the kernel and can therefore be
disregarded.)

The first application is in type $A$: our main theorem implies 
that the problem of finding a positive combinatorial formula for 
the above Gromov-Witten invariants on Grassmannians 
is a special case of the problem of providing such a formula for the 
classical structure constants on flag varieties. Both  of these open
questions have received much attention in recent years in the contexts
of quantum cohomology \cite{Ber,BCF}, 
representation theory \cite{BCM,BKMW,BMW,Tu}, 
and Schubert calculus \cite{BS,FK,S}. The present
connection becomes particularly interesting due 
to a conjecture of Knutson for the classical structure constants
\cite{Kn}.  While this conjecture is known to be
false in general, we have strong computational evidence suggesting
that the conjecture is true in the special case of two-step flag
varieties.  Thus Knutson's conjecture for two-step flags 
specializes to a conjectural quantum Littlewood-Richardson rule.

In types $B$, $C$, and $D$, our result is used to give new proofs of the
structure theorems for the quantum cohomology of Lagrangian and
orthogonal Grassmannians (denoted $LG$ and $OG$) 
which were recently obtained in 
\cite{KTlg} and \cite{KTorth}.  The strategy here is to employ
Schubert calculus on isotropic Grassmannians to compute a few key
Gromov-Witten invariants. Combining this with some basic
geometric results from the aforementioned papers, we deduce the 
quantum Pieri and Giambelli formulas for $LG$ and $OG$. 
The earlier proofs relied on an analysis of special loci on certain
Quot scheme compactifications of the moduli space of  degree $d$
rational maps to either space. In contrast, the arguments provided here
avoid the use of sophisticated moduli space constructions.

The methods of this article can also be used to prove a quantum
Pieri rule for the quantum cohomology of sub-maximal isotropic
Grassmannians. This application will be presented in \cite{BKT2}.

Our paper is organized by treating different Lie types in different
sections.  Type $A$ Grassmannians are described in Section \ref{gkn},
Lagrangian Grassmannians (type $C$) in Section \ref{lg},
and orthogonal Grassmannians
(types $B$ and $D$) in Section \ref{og}.

The authors wish to thank Michel Brion, Allen Knutson, and Laurent
Manivel for helpful discussions related to this work.

\section{The Grassmannian $G(k,n)$}
\label{gkn}

\subsection{Preliminaries}
\label{gknpre}
We review some standard definitions related to intersection theory
on Grassmannians; see \cite{F} for a detailed exposition.
Let $X = G(k,E) = G(k,n)$ be the Grassmannian
of $k$-dimensional subspaces of $E = \C^n$, and set $r = n-k$.
Given a fixed full flag $F_1 \subset F_2 \subset \cdots \subset F_n =
E$ and a partition $\lambda = (\lambda_1 \gequ \lambda_2 \gequ \cdots
\gequ \lambda_k)$ with $\lambda_1 \lequ r$ there is the Schubert variety
\[ X_\lambda(F_\bull) = \{ \,V \in X \mid \dim(V \cap
   F_{r+i-\lambda_i}) \gequ i, ~\forall \, 1 \lequ i \lequ k \,\} \,.
\]
The codimension of this variety is the weight $|\lambda| = \sum
\lambda_i$ of $\lambda$.  We let $\s_\lambda$ denote the class of
$X_\lambda(F_\bull)$ in the cohomology ring $H^*(X;\Z)$.

A rational map of degree $d$ to $X$ is a morphism
$f\colon \bP^1\to X$ such that
$$\int f_*[\bP^1]\cdot \sigma_1=d,$$ 
i.e.\ $d$ is the number of points in
$f^{-1}(X_1(F_\bull))$ when $F_\bull$ is in general position. 
All the Gromov-Witten invariants considered
in this paper are three-point and genus zero.
Given a degree $d\gequ 0$ and partitions $\lambda$, $\mu$, and $\nu$ such that
$|\lambda| + |\mu| + |\nu| = kr + nd$, we define the Gromov-Witten invariant
$\langle \s_\lambda, \s_\mu, \s_\nu \rangle_d$ to be the
number of rational maps $f\colon \bP^1\to X$ of degree $d$ such that
$f(0)\in X_\lambda(F_\bull)$, $f(1)\in X_\mu(G_\bull)$,
and $f(\infty)\in X_\nu(H_\bull)$, for given flags $F_\bull$,
$G_\bull$, and $H_\bull$ in general position.

The next proposition implies that the Gromov-Witten invariant
which counts degree $d$ maps to the Grassmannian $G(d,2d)$
through three general points, is equal to $1$.
We will say that two points $U$, $V$ of $G(d,2d)$
are in general position if the intersection $U\cap V$ of the corresponding
subspaces is trivial.

\begin{prop}
\label{3pts:typea}
Let $U$, $V$, and $W$ be three points of $\zZ=G(d,2d)$ which are
pairwise in general position.
Then there is a unique morphism $f\colon \bP^1\ra \zZ$ of degree $d$ such that
$f(0)=U$, $f(1)=V$, and $f(\infty)=W$.
\end{prop}
\begin{proof}
Let $U$, $V$, and $W$ be given, in pairwise general position.
Given a basis $\{v_1, \dots, v_d\}$ of $V$ we can construct a
morphism $f \colon \bP^1 \to \zZ$ of degree $d$ such that $f(0) = U$,
$f(1) = V$, and $f(\infty) = W$ as follows.  For $1 \lequ i \lequ d$
we let $u_i \in U$ and $w_i \in W$ be the projections of $v_i \in
\C^{2d} = U \oplus W$.  Then $f(s \cc t) = \Span \{ s u_1 + t w_1,
\dots, s u_d + t w_d \}$ satisfies the required conditions, 
where $(s \cc t)$ are the homogeneous coordinates on $\bP^1$.
Furthermore, one checks easily that the map $f$ does not
depend on the choice of  basis of $V$.

Now suppose $f \colon \bP^1 \to \zZ$ is any morphism of
degree $d$ which sends $0$, $1$, $\infty$ to $U$, $V$, $W$, respectively.
Let $S \subset \C^{2d} \otimes \cO_{\zZ}$ be the tautological bundle on $\zZ$.
The general position hypothesis implies that $f^*S = \cO(-1)^{\oplus
d}$.  Therefore we can write
\[
f(s \cc t) = \Span \{ s u_1 + t w_1,
\dots, s u_d + t w_d \}
\]
for suitable vectors $u_i, w_i \in
\C^{2d}$.  It follows that $f$ is the map constructed as above from the basis
$\{v_1,\dots,v_d\}$, where $v_i=u_i+w_i$.
\end{proof}

\medskip

Let $C$ be a non-empty subset of $X$.  Define the {\em kernel\/} of $C$
to be the intersection of all the subspaces $V \subset E$ corresponding to
points of $C$.  Similarly, the {\em span\/} of $C$ is the linear span of
these subspaces.
\[ \Ker(C) = \bigcap_{V \in C} V; \hspace{.5cm}
   \Span(C) = \sum_{V \in C} V.
\]
The kernel and span of a rational map $f\colon \bP^1\to X$
are defined to be the kernel and
span, respectively, of the image of $f$.
By considering the splitting type of the pullback of the tautological
subbundle, one proves (see \cite[Lemma 1]{Buch}) that the kernel of $f$
(resp.\ the span of $f$)
has dimension at least $k-d$ (resp.\ at most $k+d$) when $f$ has degree $d$.
In fact the main point of \cite{Buch} is to use the span of a rational map
to obtain greatly simplified proofs of type $A$ results
such as the quantum Giambelli and Pieri formulas, first proved by
Bertram \cite{Ber}.

\subsection{Main result for $G(k,n)$}
\label{gknmain}
Let $X=G(k,n)$, and fix $E=\C^n$ as before. 
Given integers $a$ and $b$, we let
$F(a,b;E)=F(a,b; n)$ denote the
two-step flag variety parametrizing pairs of subspaces $(A,B)$ with
$A \subset B \subset E$, $\dim A = a$ and $\dim B = b$.
(This variety is empty unless $0\lequ a\lequ b\lequ n$, and when 
non-empty, has dimension $(n-b)b+(b-a)a$.) 
For  any non-negative integer $d$ we set $Y_d = F(k-d, k+d; E)$.
Our main result will be used to identify Gromov-Witten invariants
on $X$ with classical intersection numbers on the flag varieties $Y_d$.

To any subvariety $\W \subset X$ we associate the following
subvariety in $Y_d$:
\begin{equation}
\label{assocsubvar}
\W^{(d)} = \{\,(A,B) \in Y_d \, \mid\, \exists~V \in \W : A \subset V
   \subset B \,\} \,.
\end{equation}
Let $F(k-d,k,k+d;E)$ denote the variety of three-step flags in $E$ of the 
indicated dimensions and 
$\pi_1 \colon F(k-d,k,k+d;E) \to X$ and $\pi_2\colon
F(k-d,k,k+d;E) \to Y_d$ the natural projections. We have $\W^{(d)} =
\pi_2(\pi_1^{-1}(\W))$; moreover, if $\W =
X_\lambda(F_\bull)$ is a Schubert variety in $X$, then
$\W^{(d)} = X_\lambda^{(d)}(F_\bull)$ is a Schubert variety
in $Y_d$ (we describe this Schubert variety in more detail after 
Corollary \ref{cor:typea}).
Since the fibers of $\pi_2$ are isomorphic to $G(d,2d)$,
the codimension of $X_\lambda^{(d)}(F_\bull)$ in $Y_d$ is
at least $|\lambda|-d^2$. Notice also that $\dim Y_d=kr+nd-3d^2$.

\begin{thm} \label{thm:typea}
Let $\lambda$, $\mu$, and $\nu$ be partitions and $d$ an integer
such that $|\lambda| + |\mu| + |\nu| = kr + nd$,
and let $F_\bull$, $G_\bull$, and $H_\bull$ be complete flags of
$E = \C^n$ in general position.  Then the
map $f \mapsto (\Ker(f),\Span(f))$ gives a bijection of the set of
rational maps $f\colon \bP^1 \to G(k,n)$ of degree $d$ satisfying
$f(0)\in X_\lambda(F_\bull)$, $f(1)\in X_\mu(G_\bull)$, and
$f(\infty)\in X_\nu(H_\bull)$, with the set of points 
in the intersection
$X^{(d)}_\lambda(F_\bull) \cap X^{(d)}_\mu(G_\bull) \cap
X^{(d)}_\nu(H_\bull)$ in $Y_d=F(k-d,k+d;n)$.
\end{thm}

\begin{proof}
Let $f\colon \bP^1\to X$ be a rational map as in the statement of the theorem.
The following dimension-counting argument shows that
$d\lequ \min(k,n-k)$, and moreover, $\dim\Ker(f)=k-d$ and $\dim\Span(f)=k+d$.
Let $a=\dim\Ker(f)$ and $b=\dim\Span(f)$.
In the two-step flag variety $Y'=F(a,b;E)$ there are 
associated Schubert varieties $X'_\lambda(F_{\bull})$, $X'_\mu(G_{\bull})$, 
and $X'_\nu(H_{\bull})$, defined as in (\ref{assocsubvar}).
If we write $e_1=k-a$ and $e_2=b-k$, then the codimension of
$X'_\lambda$ in $Y'$ is at least  $|\lambda|-e_1e_2$, 
and similar inequalities hold with $\mu$ and $\nu$ in place of $\lambda$.
Since $(\Ker(f),\Span(f))\in X'_\lambda(F_{\bull})\cap X'_\mu(G_{\bull})\cap
X'_\nu(H_{\bull})$ and the three flags are in general position, we obtain
$kr+dn-3e_1e_2\lequ \dim F(a,b;E)$, and therefore
\begin{equation}
\label{ineq1}
dn\lequ e_1a+e_2(n-b)+2e_1e_2.
\end{equation}
We know that $e_1\lequ d$ and $e_2\lequ d$, and hence that the 
right-hand side of
(\ref{ineq1}) is at most $d(a+n-b)+2e_1e_2$. It follows that 
$(e_1+e_2)^2\lequ 2d(e_1+e_2)\lequ 4e_1e_2$, and therefore
$e_1=e_2=d$.

Let $\CC$ denote the 
set of rational maps in the statement of the theorem,
and set $\II = X^{(d)}_\lambda(F_\bull) \cap
X^{(d)}_\mu(G_\bull) \cap X^{(d)}_\nu(H_\bull)$.
If $f \in \CC$ then $(\Ker(f),\Span(f)) \in \II$.
It remains to show that
for any point $(A,B) \in \II$ there is a unique map
$f \in \CC$ such that $\Ker(f)=A$ and $\Span(f)=B$.

Consider the three-step flag variety $Y''=F(k-d,k-d+1,k+d; E)$ and the 
 projection $\pi\colon Y''\ra Y_d$.  Note that $\dim Y''=\dim Y_d+2d-1$.
To each subvariety $\W\subset G(k,E)$ we associate $\W''\subset Y''$ defined 
by
\[
\W'' = \{\,(A,A',B) \in Y'' \,\mid\, \exists~V \in \W : A' \subset V
   \subset B \,\} \,.
\]
We find that the codimension of $X''_\mu(G_{\bull})$ in $Y''$ is at least 
$|\mu|-d^2+d$, and similarly for $X''_\nu(H_{\bull})$. 
Since the three flags are in general position, and
$\pi^{-1}(X^{(d)}_{\la}(F_{\bull}))$ has codimension at least $|\la|-d^2$ in $Y''$,
we have
\begin{equation}
\label{triplecap}
\pi^{-1}(X^{(d)}_{\la}(F_{\bull})) 
\cap X''_\mu(G_{\bull}) \cap X''_\nu(H_{\bull})=\emptyset,
\end{equation}
and the same is true for the other two analogous triple intersections.

Given $(A,B) \in \II$,
we let $\zZ = G(d,B/A) \subset X$ be the set of $k$-planes in $E$ between
$A$ and $B$.  Then
$X_\lambda(F_\bull) \cap \zZ$, $X_\mu(G_\bull) \cap \zZ$, and
$X_\nu(H_\bull) \cap \zZ$ are non-empty Schubert varieties in
$\zZ$. Choose three points $U$, $V$, and $W$ in $\zZ$, one from each 
intersection; the equations (\ref{triplecap}) show that these three points 
are in pairwise general position. 
Observe that any positive dimensional Schubert variety in $\zZ$ must contain
a point $U'$ which meets $U$  non-trivially, and similarly for 
$V$ and $W$. We deduce that each of
 $X_\lambda(F_\bull) \cap \zZ$, $X_\mu(G_\bull) \cap \zZ$, and
$X_\nu(H_\bull) \cap \zZ$ must be a single point.
Proposition \ref{3pts:typea} now supplies the unique $f\colon \bP^1\to X$
in $\CC$ with $\Ker(f)=A$ and $\Span(f)=B$.
\end{proof}

\begin{remark}
One can rephrase Theorem \ref{thm:typea} using
rational curves in $X$, instead of rational maps to $X$. 
Indeed, the construction in Proposition \ref{3pts:typea}
shows that every rational map $f$ that is counted 
in Theorem \ref{thm:typea} is an  
embedding of $\bP^1$ into $X$ of degree equal to the 
degree of the image of $f$. Notice also that the theorem 
implies that all of these maps have different images.
\end{remark}

It follows from Theorem \ref{thm:typea} that 
we can express any 
Gromov-Witten invariant of degree $d$ on $G(k,n)$ as
a classical intersection number on $Y_d = F(k-d,k+d;n)$.  Let
$[X_\lambda^{(d)}] \in H^*(Y_d,\Z)$ denote the cohomology class of
$X^{(d)}_\lambda(F_\bull)$.

\begin{cor}
\label{cor:typea}
Let $\lambda$, $\mu$, and $\nu$ be partitions and $d\gequ 0$  an integer
such that $|\lambda| + |\mu| + |\nu| = kr + nd$. We then have
\[
\langle \s_\lambda, \s_\mu, \s_\nu \rangle_d=
\int_{F(k-d,k+d;n)} 
[X^{(d)}_{\lambda}]\cdot [X^{(d)}_{\m}]\cdot [X^{(d)}_{\nu}].
\]
\end{cor}

\medskip

The Schubert varieties in $F(a,b;n)$ are indexed by permutations $w\in S_n$ 
with $w(i)<w(i+1)$ for all $i\notin \{a,b\}$. Given a full flag 
$F_{\bull}$ in $E$,
we define $X_w(F_{\bull})\subset F(a,b;n)$ 
as the locus of flags $A\subset B \subset E$ such that
\[
\dim(A\cap F_i) \gequ \#\{p\lequ a\ |\ w(p) > n-i\} \ \ \mathrm{and} \ \ 
\dim(B\cap F_i) \gequ \#\{p\lequ b\ |\ w(p)> n-i\}
\]
for all $i$. Recall that every partition $\la$ indexing a 
Schubert variety $X_\lambda(F_\bull) \subset G(k,n)$ 
corresponds to a Grassmannian permutation $w=w_{\la}\in S_n$, 
determined by the relations
$w(i)=\la_{k-i+1}+i$ for $i\lequ k$ and $w(i)<w(i+1)$ 
for $i>k$.
The permutation $w_{\la,d}$ associated to the modified 
Schubert variety $X^{(d)}_\la(F_\bull)$ in $F(k-d,k+d;n)$
is obtained from $w_{\la}$ by sorting the values 
$w(k-d+1),\ldots,w(k+d)$ to be in increasing order.

The proof of the next corollary gives a geometric explanation for
Yong's upper bound for the $q$-degrees in a quantum product \cite{Y}.

\begin{cor}[Yong]
\label{cor2:typea}
In the situation of Corollary \ref{cor:typea}, if any of 
$\lambda_d$, $\mu_d$, and $\nu_d$ is less than $d$, then
$\langle \s_\lambda, \s_\mu, \s_\nu \rangle_d=0$.
\end{cor}
\begin{proof}
Observe that when $\lambda_d<d$,
the codimension of
$X^{(d)}_\lambda(F_\bull)$ in $Y_d$ is strictly greater than $|\lambda|-d^2$. 
This follows by computing the length of the Weyl group element $w_{\la,d}$
associated to  $X^{(d)}_{\lambda}(F_{\bull})$. 
Therefore, when any of $\lambda_d$, $\mu_d$, or $\nu_d$ is less than $d$,
the sum of the codimensions of the three modified Schubert varieties 
which appear in the statement of Theorem \ref{thm:typea} is greater 
than the dimension of $Y_d=F(k-d,k+d;n)$.
\end{proof}

\begin{remark}
  For flag varieties $X$ other than Grassmannians, it is not in
  general true that the rational maps to $X$ counted by a
  Gromov-Witten invariant are distinguished by their kernels and spans
  (see \cite{qmonk} for definitions).  For example, if $F_\bull$,
  $G_\bull$, and $H_\bull$ are three general points in the complete
  flag variety $X = F(5)$, then there are exactly two rational maps of
  multidegree $(2,3,3,2)$ through these points, and they both have
  kernel $(0,0,0,F_4 \cap G_4 \cap H_4)$ and span
  $(F_1+G_1+H_1,\C^5,\C^5,\C^5)$.
\end{remark}

\subsection{Knutson's conjecture for two-step flag varieties}
\label{knutsonconj}

We denote by $[X_w]$ the cohomology class of the Schubert 
variety $X_w \subset F(a,b;n)$.
Each indexing permutation $w$ corresponds to a
string $J(w)$ of ``0''s, ``1''s, and ``2''s of length $n$: 
the positions of the ``0''s (resp.\
 ``1''s) in $J(w)$ are recorded by $w(1),\ldots, w(a)$ (resp.\
$w(a+1),\ldots, w(b)$).

Define a {\em puzzle\/} to be a triangle decomposed into {\em puzzle
pieces\/} of the types displayed below.  Of these pieces the fourth
and the sixth pieces come in different lengths.  The fourth piece can
have any number of ``2''s (including none) to the right of the ``0''
on the top edge and equally many to the left of the ``0'' on the
bottom edge.  Similarly the sixth piece can have an arbitrary number
of ``0''s on the top and bottom edges.

$$\pzp{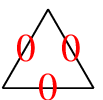} \pzs \pzp{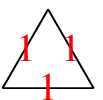} \pzs \pzp{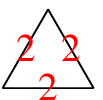} \pzs \pzp{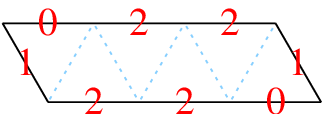}
\pzs \raisebox{-2.5mm}{\pzp{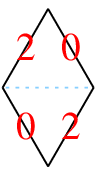}} \pzs \pzp{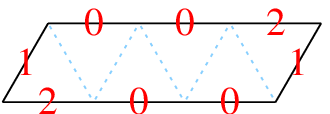}$$

\medskip

\begin{figure}
$$\includegraphics[scale=.9]{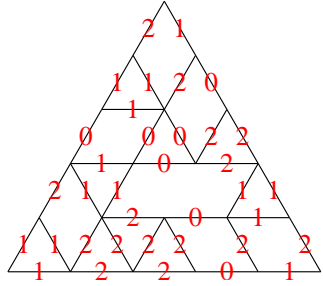}$$
\caption{An example of a puzzle}
\label{puzzleexa}
\end{figure}

A puzzle piece may be rotated but not reflected when used in a puzzle.
Furthermore, the common edges of two puzzle pieces next
to each other must have the same labels.  Figure \ref{puzzleexa} shows
an example of a puzzle.

We can now state Knutson's conjecture in the case of two-step flag
varieties.

\begin{conj}[Knutson] \label{conj:knutson}
For any three Schubert varieties $X_u$, $X_v$, and $X_w$ 
in the flag variety $F(a,b;n)$, the integral $\int_{F(a,b;n)}
[X_u]\cdot [X_v]\cdot [X_w]$ is equal to the number of puzzles such that  
$J(u)$, $J(v)$, and
$J(w)$ are the labels on the north-west, north-east, and
south sides when read in clockwise order.
\end{conj}

This conjecture has been verified for all two-step flag varieties
$F(a,b;n)$ for which $n \lequ 16$. It is a special case of a
general conjecture by Knutson for the Schubert structure constants on
all partial flag varieties, which was circulated in the fall of 1999
\cite{Kn}.  Although the Grassmannian case of this conjecture has been 
proved \cite{KTW}, Knutson quickly discovered counterexamples to the 
general conjecture.  In fact, it fails for the three-step flag variety
$F(1,3,4;5)$.


\subsection{A quantum Littlewood-Richardson rule}
\label{sec:qlr}

Let $X_\lambda(F_\bull)$ be a Schubert variety in the
Grassmannian $X = G(k,n)$ and fix a degree $d \gequ 0$. 
The $012$-string $J^d(\lambda)=J(w_{\la,d})$ associated to the 
modified 
Schubert variety $X^{(d)}_\la(F_\bull)$ in $Y_d = F(k-d,k+d;n)$
may be obtained as follows. Begin
by drawing the Young diagram of the partition $\lambda$ in the
upper-left corner of a $k \times r$ rectangle.  On the path from the
lower-left to the upper-right corner of this rectangle which follows
the border of $\lambda$ we then put a label on each step.  The
leftmost $d$ horizontal steps are labeled ``1'', while the remaining
$r-d$ horizontal steps are labeled ``2''.  The top $d$ vertical steps
are labeled ``1'', while the bottom $k-d$ vertical steps are labeled
``0''.  The string $J^d(\lambda)$ then consists of these labels in
lower-left to upper-right order.  We will write 
$X_{J^d(\lambda)}(F_\bull)$ for $X^{(d)}_\lambda(F_\bull)$.

\begin{figure}
\[
\hspace{-1.95cm}
\raisebox{-1.9cm}{\includegraphics[scale=0.9]{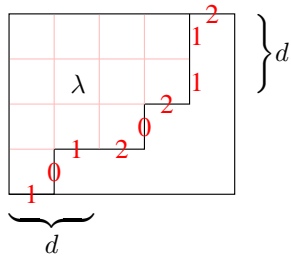}}
\hspace{0.2cm} \Biggr\}d
\hspace{-3.78cm} \underbrace{ \raisebox{-2.0cm}{\mbox{\hspace{1.13cm}}}
                                     }_{\displaystyle d}
\]
\caption{The $012$-string associated to $\lambda$ and $d$}
\label{jamsfig2}
\end{figure}

\begin{exa}
When $X = G(4,9)$, $d = 2$, and $\lambda = (4,4,3,1)$, we get
$X^{(2)}_\lambda(F_\bull) = X_{J^2(\lambda)}(F_\bull)
\subset F(2,6;9)$, where $J^2(\lambda) = 101202112$. 
This is illustrated in Figure \ref{jamsfig2}.
\end{exa}

\medskip

When Conjecture \ref{conj:knutson} is combined with Corollary
\ref{cor:typea}, we arrive at the following quantum Littlewood-Richardson rule.

\begin{conj}
\label{qlrconj}
  For partitions $\lambda,\mu,\nu$ such that $|\lambda|+|\mu|+|\nu| =
  kr + nd$ the Gromov-Witten invariant $\langle\s_\lambda,
    \s_\mu, \s_\nu\rangle_d$ is equal to the number of puzzles such
  that $J^d(\lambda)$, $J^d(\mu)$, and $J^d(\nu)$ are the labels on the
  north-west, north-east, and south sides when read in clockwise
  order.
\end{conj}

The verified cases of Conjecture~\ref{conj:knutson} imply that this
quantum Littlewood-Richardson rule is true for all Grassmannians
$G(k,n)$ for which $n \lequ 16$.  It has also been proved in some
special cases including when $d=\min(k,r)$ or when $\lambda$ has
length at most two or is a `hook' partition.  Furthermore, Conjecture
\ref{qlrconj} holds when $k \lequ 3$.  Details about this will be given
in \cite{BKT}.  See \cite{Tu} and references therein for alternative
positive formulas for some of these special cases of Gromov-Witten
invariants.

\begin{exa}
   Let $X = G(2,5)$, $d = 1$, and $Y_1 = F(1,3;5)$.  We can then compute
 \[
\langle\s_{(2,2)}, \s_{(2,1)}, \s_{(3,1)}\rangle_1 = \int_{Y_1}
  [X_{12012}] \cdot [X_{10212}] \cdot [X_{10221}] = 1.
\]
The puzzle displayed in Figure \ref{puzzleexa} is the unique puzzle with 
the strings of these Schubert classes on its sides.
\end{exa}

\section{The Lagrangian Grassmannian $LG(n,2n)$}
\label{lg}

\subsection{Preliminaries}
\label{lgprelims} In the following sections
we will mostly adopt the notational conventions of \cite{KTlg} and \cite{KTorth}. 
We begin working with the Lagrangian Grassmannian $LG=LG(n,2n)$ 
parametrizing Lagrangian (i.e.\ maximal isotropic) subspaces of 
$E=\C^{2n}$ equipped with a symplectic form. More generally, for each
$k\lequ n$, the Grassmannian $IG(k,2n)$ parametrizes $k$-dimensional
isotropic subspaces in $E$. The dimension of the complex manifold
$IG(k,2n)$ is $2k(n-k)+k(k+1)/2$.

Let $\D_n$ be the set of strict partitions $\la=(\la_1>\la_2>\cdots >\la_{\ell}>0)$ with
$\la_1\lequ n$. By convention, $\la_j=0$ for $j>\ell$. Fix an isotropic flag 
\[
F_{\bull}\ :\ 0\subset F_1\subset F_2 \subset \cdots \subset F_n \subset E
\]
of subspaces of $E$. For each $\la\in\D_n$ we have a codimension $|\la|$
Schubert variety 
$X_{\la}(F_{\bull})\subset LG$ defined as the locus of $V \in LG$ such that 
\begin{equation}
\label{lgdef}
\dim(V\cap F_{n+1-\la_i})\gequ i,\,\,i=1,\ldots,\ell(\la).
\end{equation}
Let $\s_{\la}$ be the class of $X_{\la}(F_{\bull})$ in the cohomology ring 
$H^*(LG,\Z)$.

As in Section \ref{gkn}, 
a rational map to $LG$ means a morphism $f\colon \bP^1\to LG$, and its
degree is the degree of $f_*[\bP^1]\cdot\s_1$. We also define 
the kernel  of a subset of $LG$  
as the intersection of the corresponding vector subspaces,
and the kernel  of a map $f$ to $LG$ as the kernel  
of its image (the span is the orthogonal complement of the kernel and
is not necessary here). If $f$ has degree $d$ then $\dim\Ker(f)\gequ n-d$.
The following proposition is proved in the
same way as its counterpart in type $A$.

\begin{prop}
\label{3pts:typec}
Let $U$, $V$, and $W$ be three points of $LG(d,2d)$
in pairwise general position.
Then there is a unique morphism $f\colon \bP^1\ra LG(d,2d)$ of
degree $d$ such that $f(0)=U$, $f(1)=V$, and $f(\infty)=W$.
\end{prop}

Let $h(n,d)=(n+1)(n/2+d)$. The Gromov-Witten invariant
$\langle \s_\lambda, \s_\mu, \s_\nu\rangle_d$ is defined for
$|\lambda|+|\mu|+|\nu|=h(n,d)$ and counts the number of
rational maps $f\colon\bP^1\ra LG(n,2n)$ of degree $d$ such that 
$f(0)\in X_\lambda(F_\bull)$, $f(1)\in X_\mu(G_\bull)$,
and $f(\infty)\in X_\nu(H_\bull)$, for given flags $F_\bull$,
$G_\bull$, and $H_\bull$ in general position.
The next result for the three-point Gromov-Witten invariants of 
degree $1$ on $LG$ was proved in \cite[\S 3]{KTlg} via a
geometric correspondence between lines on $LG(n,2n)$ and
points in $LG(n+1,2n+2)$:

\begin{prop}
\label{linenumbers}
For $\la$, $\m$, $\n\in \D_n$ we have 
\[
\langle \sigma_{\la},\sigma_\m,\sigma_\n\rangle_1=
\frac{1}{2}\int_{LG(n+1,2n+2)}[X^+_\la]\cdot[X^+_\m]\cdot[X^+_\n],
\]
where $X_{\la}^+$, $X^+_\m$, $X^+_\n$ 
 denote Schubert varieties in $LG(n+1,2n+2)$.
\end{prop}

\subsection{Main result for $LG(n,2n)$}
\label{mainLG}

For any integer $d$ with $0\lequ d\lequ n$, 
let $Y_d$ denote the isotropic Grassmannian 
$IG(n-d,2n)$. Given a subvariety $\W\subset LG$, we define
$\W^{(d)}\subset Y_d$ by the prescription
\begin{equation}
\label{cwdef}
\W^{(d)}=\{\,A\in Y_d\, \mid\,
\exists~V \in \W : A \subset V\,\}\,.
\end{equation}
When $\W=X_\lambda(F_\bull)$ is a Schubert variety in $LG$,
we thus obtain a Schubert variety $X^{(d)}_\lambda(F_\bull)$
in $Y_d$, of codimension at least $|\la|-d(d+1)/2$.
Note that $X^{(1)}_\lambda$ 
and $X^{(2)}_\lambda$ are the varieties which were denoted
$\X_{\lambda}'$ and $\X_{\lambda}''$, respectively, in \cite[\S 3.1]{KTlg}.

\begin{thm}
\label{thm:typec}
Let $d$ be an integer and $\lambda$, $\mu$, $\nu\in\D_n$ be
such that $|\lambda|+|\mu|+|\nu|=h(n,d)$,
and let $F_\bull$, $G_\bull$, and $H_\bull$ be complete isotropic flags
of $E=\C^{2n}$ in general position.
Then the map $f\mapsto\Ker(f)$ gives a bijection of the set of
rational maps $f\colon \bP^1 \to LG(n,2n)$ of degree $d$ satisfying
$f(0)\in X_\lambda(F_\bull)$, $f(1)\in X_\mu(G_\bull)$, and
$f(\infty)\in X_\nu(H_\bull)$, with the set of points 
in the intersection
$X^{(d)}_\lambda(F_\bull) \cap X^{(d)}_\mu(G_\bull) \cap
X^{(d)}_\nu(H_\bull)$ in $Y_d=IG(n-d,2n)$.
\end{thm}
\begin{proof}
We proceed as in the proof of Theorem \ref{thm:typea}.
Given a map $f$ as in the statement of the theorem, set
$e=n-\dim\Ker(f)$; we know that $e\lequ d$. The codimension of
$X^{(e)}_\lambda(F_{\bull})$ in $Y_e$ is at least $|\lambda|-e(e+1)/2$,
and similarly for $X^{(e)}_{\mu}(G_{\bull})$ and $X_{\nu}^{(e)}(H_{\bull})$.
Since $\Ker(f)\in X^{(e)}_{\la}(F_{\bull})\cap X^{(e)}_{\mu}(G_{\bull})\cap
X^{(e)}_{\nu} (H_{\bull})$, a dimension count shows that 
\[
h(n,d)-3e(e+1)/2 \lequ \dim IG(n-e,2n),
\]
which simplifies to $d\lequ e$,
and hence $\dim\Ker(f)=n-d$. Therefore, $\Ker(f)\in\II$, where
$\II$ denotes the intersection in $Y_d$ from the
statement of the theorem. The argument is completed
as in type $A$, except that the type $C$ analogue of the 
variety $Y''$  is the isotropic partial flag space $IG(n-d,n-d+1;2n)$, 
and we invoke Proposition \ref{3pts:typec} 
instead of Proposition \ref{3pts:typea}. 
\end{proof}

Exactly as in type $A$, Theorem \ref{thm:typec} identifies 
Gromov-Witten invariants on
$LG$ with classical structure constants on isotropic Grassmannians. 

\begin{cor}
\label{cor:typec}
Let $d\gequ 0$ and choose
$\la$, $\m$, $\n\in\D_n$ with 
$|\la|+|\m|+|\n|=h(n,d)$. We then have
\[
\langle \s_\la, \s_\m, \s_\n \rangle_d=
\int_{IG(n-d,2n)} [X^{(d)}_{\la}]\cdot [X^{(d)}_{\m}]\cdot [X^{(d)}_{\nu}].
\] 
\end{cor}

\medskip

Each $\la\in \D_n$ corresponds to a maximal Grassmannian element
\[
w_{\la}=(\ov{\la}_1,\ldots,\ov{\la}_{\ell},\la'_{n-\ell},\ldots,\la'_1)
\]
of the Weyl group $W_n$ for the root system $C_n$ (see \cite[Sec.\ 2]{KTlg}).
When the length $\ell$ of $\la$ is at least  $d$, the Weyl group element
$w_{\la,d}$ in $W_n$ associated to $X^{(d)}_\la(F_\bull)$ 
is obtained from $w_{\la}$ by removing the bars from the first $d$ entries and
sorting them to be in increasing order.
In this case, the codimension of $X^{(d)}_\la(F_{\bull})$ 
in $IG(n-d,2n)$ is exactly $|\la|-d(d+1)/2$. 
 
The vanishing statement in next corollary was shown 
in \cite[Prop.\ 10]{KTlg} using different methods; we argue here as in the proof
of Corollary \ref{cor2:typea}.

\begin{cor}
\label{cor2:typec}
In the situation of Corollary \ref{cor:typec}, if
any of $\la$, $\m$, and $\n$ has length less than $d$, then
 $\langle \s_\lambda, \s_\mu, \s_\nu \rangle_d=0$.
\end{cor}

\subsection{Quantum cohomology of $LG(n,2n)$}
\label{lgapp}

In this section, we use Corollaries \ref{cor:typec} and \ref{cor2:typec}
to give new proofs of the quantum Pieri rule \cite[Prop.\ 8]{KTlg} and
quantum Giambelli formula \cite[Thm.\ 1]{KTlg} for $LG$. These are key
ingredients used to describe the multiplicative structure of the quantum
cohomology ring $QH^*(LG(n,2n))$. Our proofs will require the crucial fact
that the product defined in (\ref{qprod}) is associative; a reference for
this basic result is the expository paper \cite{FP}.


The (small) quantum cohomology ring $QH^*(LG)$ is a $\Z[q]$-algebra
which is isomorphic to $H^*(LG,\Z)\otimes_{\Z}\Z[q]$ as a module over
$\Z[q]$. Here $q$ is a formal variable of degree $n+1$. The ring
structure on $QH^*(LG)$ is determined by the relation
\begin{equation}
\label{qprod}
\s_{\la}\cdot \s_{\m} = 
\sum \langle \s_{\la}, \s_{\m}, \s_{\n'} \rangle_d\,\s_{\n}\, q^d,
\end{equation}
the sum over $d\gequ 0$ and strict partitions $\n$ with
$|\n|=|\la|+|\m|-d(n+1)$. Here $\n'$ denotes the
dual partition of $\n$, defined so that the classes $\s_{\n}$ and
$\s_{\n'}$ are Poincar\'e dual to each other in $H^*(LG,\Z)$.

We recall some further terminology involving partitions and Young diagrams.
A skew diagram $\alpha$ is a horizontal strip if it has at most
one box in each column; we define the connected components of $\alpha$
by agreeing that two boxes in $\alpha$ are connected if they share a
vertex or an edge. If $\mu\subset \la$, we let $N(\mu,\la)$ denote the
number of connected components of $\la/\mu$ which do not meet the
first column, and we let $N'(\m,\la)$ be the one less than the total
number of connected components of $\la/\m$. We use the notation
$\la\ssm\m$ to denote the partition whose parts are the parts of $\la$
which are not parts of $\mu$.

\begin{thm}[Quantum Pieri rule for $LG$] 
\label{qupieri}
For any $\lambda\in\D_n$ and $k\gequ 0$ we have
\begin{equation}
\label{quantumpieri}
\s_{\lambda}\cdot\s_k=\sum_{\m} 2^{N(\lambda,\m)}\s_{\m}+\sum_{\n}
2^{N'(\nu,\lambda)} \s_{\nu} \, q
\end{equation}
in $QH^*(LG(n,2n))$, where the first sum is over all strict partitions
$\m\supset\lambda$ with $|\m|=|\lambda|+k$ such that $\m/\lambda$ is a
horizontal strip, while the second is over all strict $\n$ contained
in $\lambda$ with $|\n|=|\lambda|+k-n-1$ such that $\lambda/\n$ is a
horizontal strip.
\end{thm}

\begin{proof}
The first sum in (\ref{quantumpieri}) is just a restatement of the
classical Pieri rule of Hiller and Boe \cite{HB}.  It follows from
Proposition \ref{linenumbers} that the coefficient of $q$ in
(\ref{quantumpieri}) is given by the classical Pieri rule for
$LG(n+1,2n+2)$; one checks easily that the multiplicities agree with
those in the theorem (see \cite[Prop.\ 7]{KTlg}). Finally, Corollary
\ref{cor2:typec} shows that for all $d\gequ 2$, the coefficient of
$q^d$ in $\s_{\lambda}\cdot \s_k$ is zero.
\end{proof}

\begin{thm}[Ring presentation and quantum Giambelli]
\label{lggiambthm} 
The ring $QH^*(LG)$ is presented as a quotient of the polynomial
ring $\Z[\s_1,\ldots,\s_n,q]$ by the relations 
\begin{equation}
\label{lgqrel}
\s_i^2+2\sum_{k=1}^{n-i}(-1)^k\s_{i+k}\s_{i-k}=(-1)^{n-i}\s_{2i-n-1}\,q
\end{equation}
for $1\lequ i\lequ n$, 
where it is understood that $\s_0=1$ and $\s_j=0$ for $j<0$.
The Schubert class $\s_{\lambda}$
in this presentation is given by the Giambelli formulas
\begin{equation}
\label{lgqgiam1}
\s_{i,j}=
\s_i\s_j+2\sum_{k=1}^{n-i}(-1)^k\s_{i+k}\s_{j-k}+(-1)^{n+1-i}\s_{i+j-n-1}\,q
\end{equation}
for $i>j>0$, and 
\begin{equation}
\label{lgqgiam2}
\s_{\lambda}=\text{\em Pfaffian}[\s_{\lambda_i,\lambda_j}]_{1\lequ i<j\lequ r},
\end{equation}
where multiplication in $QH^*(LG(n,2n))$ is employed
throughout  and $r$ is the smallest even integer such that $r \gequ \ell(\la)$.
\end{thm}

\begin{proof}
The quantum relations (\ref{lgqrel}) and the two-condition quantum
Giambelli formula (\ref{lgqgiam1}) are easily deduced from Proposition
\ref{linenumbers} (following \cite[\S 3.3]{KTlg}).  We give a new
proof of the Pfaffian formula (\ref{lgqgiam2}).

Suppose $\lambda\in \D_n$  with $\ell(\lambda)\gequ 3$.  
To prove (\ref{lgqgiam2}), it suffices to establish the
following identity in $QH^*(LG(n,2n))$:
\begin{equation}
\label{qgeqtwo}
\sigma_\lambda=\sum_{j=1}^{r-1} (-1)^{j-1} \sigma_{\lambda_j,\lambda_r} 
\sigma_{\lambda\smallsetminus\{\lambda_j,\lambda_r\}}.
\end{equation}
The classical component of (\ref{qgeqtwo}) is the
classical Pfaffian Giambelli identity for $LG$, as given by
Pragacz \cite{P}. Hence, (\ref{qgeqtwo}) is equivalent to
\begin{equation}
\label{qgeqthree}
0=\sum_{j=1}^{r-1} (-1)^{j-1} \langle \sigma_{\lambda_j,\lambda_r},
\sigma_{\lambda\smallsetminus\{\lambda_j,\lambda_r\}},
\sigma_\m \rangle_d
\end{equation}
for every $d\gequ 1$ and $\m\in \D_n$ such that $|\lambda|+|\m|=h(n,d)$.
By Corollary \ref{cor:typec}, the vanishing in (\ref{qgeqthree}) will
follow from the relations
\begin{equation}
\label{idone}
\sum_{j=1}^{r-1}(-1)^{j-1}\,[X^{(1)}_{\lambda_j,\lambda_r}]\,
[X^{(1)}_{\lambda\ssm\{\lambda_j,\lambda_r\}}]=0
\end{equation}
and 
\begin{equation}
\label{idtwo}
\sum_{j=1}^{r-1} (-1)^{j-1}\,[X^{(2)}_{\lambda_j,\lambda_{r}}]\,
[X^{(2)}_{\lambda\ssm\{\lambda_j,\lambda_{r}\}}]=0
\end{equation}
in $H^*(IG(n-1,2n),\Z)$ and $H^*(IG(n-2,2n),\Z)$, respectively.

The cohomological formulas (\ref{idone}) and (\ref{idtwo}) were 
both proved in  \cite[Cor.\ 1]{KTlg} by an algebraic argument which 
established more: namely, that they hold  
on the level of symplectic Schubert polynomials. The validity of precisely 
these two identities was suggested by the analysis of the Lagrangian Quot 
scheme in \cite{KTlg}. We refer to op.\ cit.\ for more details.
\end{proof}

\begin{remark}
We have given independent proofs of the quantum Pieri formula and the
quantum Giambelli formula. However, it suffices to prove either of these two
results; the other follows formally.
In \cite{KTlg}, the quantum Giambelli formula is proved using intersection
theory, and the quantum Pieri formula is deduced by an algebraic argument.
For the reverse implication, one may argue as follows. The first relation
(\ref{lgqgiam1}) is derived easily from (\ref{quantumpieri}), so we prove the 
Pfaffian formula (\ref{lgqgiam2}).  We employ the
algebra of $\wt{Q}$-polynomials \cite{PR}.  Each $\wt{Q}$-polynomial  is 
indexed by a (not necessarily strict) partition $\lambda$,
and there are integers $e(\lambda,\m;\,\n)$ such that
\[
\wt{Q}_{\lambda}\,\wt{Q}_{\mu}=\sum e(\lambda,\m;\,\n)\,\wt{Q}_{\n}.
\]
One now shows by induction on the length $\ell(\mu)$ that
\begin{equation}
\label{inductuptoqgiambelli}
\s_\lambda\cdot \s_\m=\sum 2^{-d} 
e(\lambda,\m;\,((n+1)^d,\n))\,\s_\n\, q^d,
\end{equation}
using the quantum Pieri rule and the Pieri
formula for $\wt{Q}$-polynomials (\cite[Prop. 4.9]{PR}) as the base case.
The quantum Giambelli formula (\ref{lgqgiam2}) is recovered by
combining (\ref{inductuptoqgiambelli}) with the
Pfaffian relations
\[
\wt{Q}_{\lambda}=
\sum_{j=1}^{r-1}(-1)^{j-1}\wt{Q}_{\lambda_j,\lambda_r}
\cdot \wt{Q}_{\lambda\ssm\{\lambda_j,\lambda_r\}}.
\]
\end{remark}

\section{The orthogonal Grassmannian $OG(n+1,2n+2)$}
\label{og}

\subsection{Preliminaries}
\label{ogpre}
In this section we give an analogous treatment of the even orthogonal
Grassmannian $OG=OG(n+1,2n+2)=SO_{2n+2}/P_{n+1}$. This variety
parametrizes (one component of) the locus of maximal isotropic subspaces
of a $(2n+2)$-dimensional vector space $E$, equipped with a nondegenerate
symmetric form. Recall that $OG$ is isomorphic to the odd orthogonal 
Grassmannian $SO_{2n+1}/P_n$, hence the results of this section
apply to the latter variety as well.
The Schubert varieties $X_{\lambda}(F_{\bull})$ in $OG$ are again
parametrized by partitions $\lambda\in\D_n$ and defined by the equations
(\ref{lgdef}), with respect to a complete isotropic flag $F_{\bull}$ in $E$. 
Let $\tau_{\lambda}$ be the cohomology class of $X_{\la}(F_{\bull})$; the
set $\{\tau_{\la}\, |\, \la\in\D_n\}$ is a $\Z$-basis of $H^*(OG,\Z)$.
For each $k\gequ 0$, let $OG(k,2n+2)$ denote the orthogonal Grassmannian 
parametrizing $k$-dimensional isotropic subspaces of $E$.

As before, the degree of a morphism $f\colon\bP^1\ra OG$ 
is defined as the degree of $f_*[\bP^1]\cdot \tau_1$.
We say that two points $U$, $V$ of $OG(e,2e)$ are in general position if the
corresponding subspaces have trivial intersection; the parity
condition $\dim(U\cap V) \equiv e \, ($mod $2)$ shows that
this can happen only if
$e$ is even. Proposition \ref{3pts:typec} has the following 
orthogonal analogue.

\begin{prop}
\label{3pts:typed}
Let $U$, $V$, and $W$ be three points of $OG(e,2e)$
in pairwise general position (so, in particular, $e=2d$ is even).
Then there is a unique morphism $f\colon \bP^1\ra OG(e,2e)$ of
degree $d$ such that $f(0)=U$, $f(1)=V$, and $f(\infty)=W$.
\end{prop}

An important difference between this and the previous section is that
the natural embedding of $OG(n+1,2n+2)$ 
into the type $A$ Grassmannian $G(n+1,2n+2)$
is degree-doubling. This means that for any degree $d$ map 
$f\colon\bP^1\ra OG$, the  pullback of the tautological quotient bundle 
over $OG$ has degree $2d$.
Moreover, the dimension of $\Ker(f)$ is bounded below by $n+1-2d$.

\subsection{Main result for $OG(n+1,2n+2)$}
\label{ogmain}

 Let $Y_d$ denote the orthogonal Grassmannian $OG(n+1-2d,2n+2)$.
For a given $d$ with $0\lequ d\lequ (n+1)/2$ and subvariety 
$\W\subset OG(n+1,2n+2)$,  define $\W^{(d)}\subset  Y_d$ by the  
same prescription (\ref{cwdef}) as in type $C$. For each Schubert
variety $X_{\la}=X_{\la}(F_{\bull})$ in $OG$, we thus get a Schubert variety
$X_{\la}^{(d)}=X_{\la}^{(d)}(F_{\bull})$ in $Y_d$.  We note that for $d=1$,
$X_{\la}^{(1)}$ is the variety which was denoted by $\Y_{\la}$ in 
\cite[\S 3.1]{KTorth}.

 Let $h'(n,d)=n(n+1)/2+2nd$.
The Gromov-Witten invariant $\langle \tau_\lambda,\tau_\m,\tau_\n\rangle_d$
is zero unless the dimension condition
$|\lambda|+|\m|+|\n|=h'(n,d)$ is satisfied, and this occurs only for
$d\lequ (n+1)/2$. 
We now have the following orthogonal analogue of Theorem \ref{thm:typec},
which is used to
identify Gromov-Witten invariants on $OG(n+1,2n+2)$ with classical
structure constants on non-maximal orthogonal Grassmannians.

\begin{thm}
\label{thm:typed}
Let $d$ be an integer and $\lambda$, $\mu$, $\nu\in\D_n$ be
such that $|\lambda|+|\mu|+|\nu|=h'(n,d)$,
and let $F_\bull$, $G_\bull$, and $H_\bull$ be complete isotropic flags
of $E=\C^{2n+2}$ in general position.
Then the map $f\mapsto\Ker(f)$ gives a bijection of the set of
rational maps $f\colon \bP^1 \to OG(n+1,2n+2)$ of degree $d$ satisfying
$f(0)\in X_\lambda(F_\bull)$, $f(1)\in X_\mu(G_\bull)$, and
$f(\infty)\in X_\nu(H_\bull)$, with the set of points 
in the intersection
$X^{(d)}_\lambda(F_\bull) \cap X^{(d)}_\mu(G_\bull) \cap
X^{(d)}_\nu(H_\bull)$ in $Y_d=OG(n+1-2d,2n+2)$.
\end{thm}


\begin{cor}
\label{cor:typed}
Let $d\gequ 0$ and choose
$\la$, $\m$, $\n\in\D_n$ with 
$|\la|+|\m|+|\n|=h'(n,d)$. We then have
\[
\langle \tau_\lambda, \tau_\m, \tau_\n \rangle_d=
 \int_{OG(n+1-2d,2n+2)} [X^{(d)}_{\la}]\cdot 
[X^{(d)}_{\m}]\cdot [X^{(d)}_{\nu}].
\] 
\end{cor}


\begin{cor}
\label{cor2:typed}
In the situation of Corollary \ref{cor:typed}, 
if any of $\lambda$, $\m$, and $\n$ has length less than $2d-1$, then
the Gromov-Witten invariant $\langle \tau_\lambda, \tau_\m, \tau_\n \rangle_d$
vanishes.
\end{cor}

The proofs of these results are similar to those
of Theorem \ref{thm:typec} and its corollaries.
This time the auxiliary space $Y''$ is the isotropic 
partial flag variety $OG(n+1-2d,n+3-2d;2n+2)$, and we invoke
Proposition \ref{3pts:typed}. Note that a stronger version of the vanishing 
statement in Corollary \ref{cor2:typed} was established in 
\cite[Thm.\ 6]{KTorth}.

\subsection{Quantum cohomology of $OG(n+1,2n+2)$}

The quantum cohomology ring of $OG(n+1,2n+2)$ is defined in a similar way 
to that of $LG$. A big difference here is that the degree of $q$ is $2n$. 
We now give independent proofs of the basic
structure theorems \cite[Thm.\ 1]{KTorth}
and \cite[Cor.\ 5]{KTorth} regarding the 
quantum cohomology of $OG$. 

\medskip

\begin{thm}[Quantum Pieri rule for $OG$] 
\label{ogqupieri}
For any $\lambda\in\D_n$ and $k\gequ 0$ we have
\begin{equation}
\label{ogqupieriformula}
\tau_{\lambda}\cdot\tau_k=\sum_{\m} 2^{N'(\lambda,\m)}\tau_{\m}+\
\sum_{\nu}2^{N'(\lambda,\nu)} \tau_{\nu\ssm (n,n)} \, q,
\end{equation}
where the first sum is over all strict $\mu\supset\lambda$ with
$|\mu|=|\lambda|+k$ such that $\mu/\lambda$ is a horizontal strip, 
and the second sum is over all partitions $\nu=(n,n,\ov{\nu})$ with 
$\ov{\nu}$ strict such that $\nu\supset \la$, 
$|\nu|=|\lambda|+k$, and $\nu/\lambda$ is a horizontal strip.
\end{thm}

\begin{proof}
The argument here differs from that in the proof of Theorem \ref{qupieri}.
We first show that $\tau_\lambda\cdot\tau_k$ is classical whenever
$\lambda_1<n$; in other words, if the first row of $\lambda$ is not full, then
multiplying $\tau_{\la}$ by a special Schubert class carries no quantum 
correction. Observe that Corollary \ref{cor2:typed} implies 
$\langle \tau_\lambda, \tau_k, \tau_{\mu}\rangle_d=0$  for $d>1$.

Suppose now that some Gromov-Witten invariant
$\langle \tau_\lambda, \tau_k, \tau_{\mu}\rangle_1$ is non-zero.
It follows that for some complete isotropic flags $F_{\bull}$, $G_{\bull}$,
and $H_{\bull}$, there are some (only finitely many)
lines on $OG$ incident to  
$X_{\la}(F_{\bull})$, $X_k(G_{\bull})$, and $X_{\mu}(H_{\bull})$. 
One sees easily that given a point $U\in X_\lambda(F_{\bull})$, 
with $\lambda_1<n$, there exists a  flag $F'_{\bull}$ such that 
$U\in X_{\rho_{n-1}}(F'_{\bull})\subset
X_{\la}(F_{\bull})$, where $\rho_{n-1}=(n-1,\ldots,2,1)$.
Moreover, for any $V\in X_k(G_{\bull})$, 
there exists a flag $G'_{\bull}$ such that 
$V\in X_n(G'_{\bull})\subset X_k(G_{\bull})$.

We claim now that whenever there is one line on $OG$ incident
to $X_{\rho_{n-1}}(F_{\bull})$, $X_n(G_{\bull})$, and the point $W$,
there are infinitely many lines satisfying these incidence conditions.
In light of this, it is impossible to have
$\langle \tau_\lambda, \tau_k, \tau_\m\rangle_1\ne 0$.
To prove the claim, we recall that lines on $OG$  
are in bijective correspondence with points in
$OG(n-1,2n+2)$, and we translate the given incidence
conditions to conditions on the corresponding
$(n-1)$-dimensional isotropic space $T$:

(i) the line meets $X_n(G_{\bull})$ if and only if
$T\subset G_1^{\perp}$;

(ii) the line meets $X_{\rho_{n-1}}(F_{\bull})$ if and only if
$\dim(T\cap\wt{F}_{n+1})\gequ n-2$, where $\wt{F}_{n+1}$ is 
the unique maximal isotropic subspace containing $F_n$ other 
than $F_{n+1}$;

(iii) the line meets the point $W$ if and only if $T\subset W$.

\medskip

Assume that $T$ satisfies (i)--(iii). Then
$$\dim(G_1^{\perp}\cap \wt{F}_{n+1}\cap W)\gequ n-2.$$
Consider a fixed subspace 
$S\subset G_1^{\perp}\cap \wt{F}_{n+1}\cap W$, with $\dim S=n-2$.
Now any $T'$ of dimension $n-1$ containing $S$ and contained in
$G_1^\perp\cap W$ satisfies (i)--(iii) as well, and there are
infinitely many such $T'$.

Next, consider any product $\tau_\lambda\cdot\tau_k$.
The first sum in (\ref{ogqupieriformula}) agrees with the classical
Pieri formula, so we focus on the quantum term.
When $\lambda_1<n$, the quantum term (\ref{ogqupieriformula}) correctly
vanishes, so we suppose $\lambda_1=n$.
We write
$\lambda\ssm n=(\lambda_2,\lambda_3,\ldots)$; we then have the equation
$$\tau_\lambda=\tau_n\cdot\tau_{\lambda\ssm n}$$ 
in $QH^*(OG)$.

Suppose, first, that $k=n$. Note that the formula
$\tau_n^2=q$
holds in $QH^*(OG)$; this follows from the easy enumerative fact that
there is a unique line in $OG$ through a given point and incident to  
 $X_n(F_{\bull})$ and $X_n(G_{\bull})$, for general flags $F_{\bull}$ and
$G_{\bull}$. Therefore, in this case, we have 
$\tau_k\cdot\tau_\lambda=\tau_n^2\,\tau_{\lambda\ssm n}=
q\,\tau_{\lambda\ssm n}$,
and the quantum Pieri formula is verified.
If $k<n$, then we write
$$\tau_k\cdot\tau_\lambda=\tau_n\cdot(\tau_k\cdot\tau_{\lambda\ssm n}).$$
The product in parentheses receives no quantum correction, and hence is
given by the classical Pieri formula.
As the quantum Pieri formula for multiplication by $\tau_n$ has already
been established, it remains only to verify that the result
agrees with (\ref{ogqupieriformula}), and this is easily checked.
\end{proof}


\begin{thm}[Ring presentation and quantum Giambelli]
\label{oggiambthm} 
The ring $QH^*(OG)$ is presented as a quotient of the polynomial
ring $\Z[\tau_1,\ldots,\tau_n,q]$ modulo the relations 
\[
\tau_i^2+2\sum_{k=1}^{i-1}(-1)^k\tau_{i+k}\tau_{i-k}+(-1)^i\tau_{2i}=0
\]
for all $i<n$, together with the quantum relation
\[
\tau_n^2=q
\]
{\em(}it is understood that $\tau_j=0$ for $j>n${\em)}.
The Schubert class $\tau_{\lambda}$
in this presentation is given by the Giambelli formulas
\begin{equation}
\label{ogqgiam1}
\tau_{i,j}=
\tau_i\tau_j+2\sum_{k=1}^{j-1}(-1)^k\tau_{i+k}\tau_{j-k}+(-1)^j\tau_{i+j}
\end{equation}
for $i>j>0$, and 
\begin{equation}
\label{ogqgiam2}
\tau_{\lambda}=\text{\em Pfaffian}
[\tau_{\lambda_i,\lambda_j}]_{1\lequ i<j\lequ r},
\end{equation}
where quantum multiplication is employed throughout 
and $r$ is the smallest even integer such that $r \gequ \ell(\la)$.
\end{thm}

\begin{proof}
That the ring presentation is as claimed follows from the argument in
\cite{ST} (see also \cite[\S 3.3]{KTorth}).
Formula (\ref{ogqgiam1}) is true classically and holds without
any quantum correction for degree reasons, since $\deg(q)=2n$.
Formula (\ref{ogqgiam2}) is equivalent to the Pfaffian Laplace-type
expansion
\begin{equation}
\label{quantumgiambellieqtwo}
\tau_\lambda=\sum_{j=1}^{r-1} (-1)^{j-1} \tau_{\lambda_j,\lambda_r}
\tau_{\lambda\smallsetminus\{\lambda_j,\lambda_r\}}.
\end{equation}
Observe that (\ref{quantumgiambellieqtwo}) will follow from
Corollary \ref{cor:typed} and the formula 
\begin{equation}
\label{oequ}
\sum_{j=1}^{r-1}(-1)^{j-1}[X^{(1)}_{\lambda_j,\lambda_r}]\, 
[X^{(1)}_{\lambda\ssm\{\lambda_j,\lambda_r\}}] = 0
\end{equation}
in $H^*(OG(n-1,2n+2), \Z)$. Finally, (\ref{oequ}) was established
in \cite[Cor.\ 1]{KTorth} by an algebraic argument, whose  
motivation was similar to that in the Lagrangian case.
\end{proof}

At the end of \S \ref{lgapp} we noted that it suffices to prove either
the quantum Giambelli or quantum Pieri formula; the other then follows
using the algebra of symmetric polynomials. The same remark applies to
$OG$, and the details are similar (one of the two implications was worked 
out in \cite[Sec.\ 6]{KTorth}).

\end{document}